\documentclass[10pt]{article}
\newcommand{\Keywords}[1]{\par\noindent
{\small{\em Keywords\/}: #1}}

\usepackage{verbatim}
\usepackage{amssymb,xcolor}
\usepackage{amsmath}
\usepackage{graphicx}
\newtheorem{thm}{Theorem}

\newtheorem{lem}[thm]{Lemma}
\newtheorem{prop}[thm]{Proposition}

\newenvironment{proof}[1][Proof]{\textbf{#1.} }{\ \rule{0.5em}{0.5em}}

\def\bt{\hbox{$\bf \cdot$}}

\let\eps=\varepsilon

\let\D=\Delta

\let\s=\sigma

\let\t=\theta
\let\l=\lambda

\let\ls=\leqslant
\let\gs=\geqslant
\let\le=\ls
\let\ge=\gs

\let\cal=\mathcal
\def\dvit{\colon\ }
\def\dr{\hbox{\rm d}}

\def\E{\hbox{\bf E}}

\newcommand{\fy}{\varphi}
\newcommand{\tend}[1]{\longrightarrow}
\newcommand{\tends}[1]{\xrightarrow[#1]{}}
\newcommand{\tendsd}{\xrightarrow{\ d\ }}
\newcommand{\abs}[1]{\lvert #1\rvert}

\title{Estimation of parameters of SDE driven by fractional Brownian motion
with polynomial drift}

\date{}
\author{K. Kubilius *$^1$, V. Skorniakov $^2$, D. Melichov $^3$}

\begin{document}

\maketitle

\vskip-0.5cm

\centerline{\small $^1$ Vilnius University, Institute of Mathematics and Informatics, Akademijos 4,}
\centerline{\small LT-08663 Vilnius, Lithuania}

\centerline{\small $^2$ Vilnius University, Faculty of Mathematics and
Informatics, Naugarduko 24,} \centerline{\small LT-03225, Vilnius, Lithuania}

\centerline{\small $^3$ Vilnius Gediminas Technical University, Saul\.etekio al.
11,} \centerline{\small LT-10223, Vilnius, Lithuania}

\let\oldthefootnote\thefootnote
\renewcommand{\thefootnote}{\fnsymbol{footnote}}
\footnotetext[1]{This research was funded by a grant (No.
MIP-048/2014) from the Research Council of Lithuania.}

\abstract{Strongly consistent and asymptotically normal estimators
of the Hurst index and volatility parameters of solutions of
stochastic differential equations with polynomial drift are
proposed. The estimators are based on discrete observations of the
underlying processes.

\bigskip
 }

\section{Introduction}

Many applications use processes that are described by stochastic
differential equations (SDE). Recently, much attention  has been
paid to SDEs driven by the fractional Brownian motion (fBm) and to
the problems of statistical estimation of model parameters.
Statistical aspects of the models including fBm were studied in many
articles. Especially much attention was paid to the estimations of
the parameter of drift. We concentrate on estimation of the Hurst
parameter and volatility of SDE with polynomial drift
\begin{equation}\label{e:diflygt}
X_t=\xi+\int_0^t\big(a X_s^m+b X_s\big) ds + c\int_0^t X_s\, d
B^H_s,\quad \xi\gs0,\quad  t\in[0,T],
\end{equation}
where $1/2<H<1$, $m\in \mathbb{N}$, $m\gs 2$, and $B^H$ is a fBm with Hurst
index $1/2<H<1$. Almost all sample paths of $B^H$ have bounded
$p$-variation for each $p>1/H$ on $[0,T]$ for every $T>0$ (a more
detailed definition of an fBm will be formulated in the next
section). The second integral in (\ref{e:diflygt}) is the pathwise
Riemann-Stieltjes integral with respect to the process having finite
$p$-variation.

It is well known that equation (\ref{e:diflygt}) for standard
Brownian motion, i. e. for $H=1/2$,  is nonlinear reducible SDE with
polynomial drift of degree $m$. Under the substitution
$h(x)=x^{1-m}$ equation (\ref{e:diflygt}) reduces to a linear SDE
with multiplicative noise. Some famous equations can be obtained as
particular cases of this equation. For example,
we can obtain the Black-Scholes model, the Verhulst and the
Landau-Ginzburg equations.

Using an approximation approach in $L^2$ Dung \cite{dung} obtained
an explicit solution of the equation (\ref{e:diflygt}) for $H>1/2$. His
approximation is based on the fact that fBm of Liouville form can be
uniformly approximated by semimartingales.

Many authors (see \cite{glad}, \cite{IL}, \cite{Cohen-98},
\cite{Coeurjolly-01}, \cite{Begyn-06}, \cite{begyn2},  \cite{rn1},
\cite{ma}, \cite{kub}) considered an almost sure convergence and an
asymptotic normality of the generalized quadratic variations
associated to the filter $a$ (see \cite{IL}) of a wide class of
processes with Gaussian increments. Conditions  for this type of
convergence are expressed in terms of covariances of a Gaussian
process and depend on some parameter $\gamma\in(0,2)$. These results
make it possible to obtain an  estimator of the parameter $\gamma$
and to study its asymptotic properties. If the process under
consideration is the fBm then this parameter is the Hurst index
$0<H<1$.

There are a few results concerning estimation of the Hurst index if
the observed process is described by SDE driven by fBm.  For the
first time, to our knowledge,  an estimator of the Hurst parameter
of a pathwise solution of a linear SDE driven by a fBm and its
asymptotical behaviour were considered in \cite{Berzin-08}.  A more general
situation was considered in \cite{kubmel1} using
a different approach. Estimators for the Hurst parameter were constructed
making use of the first and second order quadratic variations of the
observed values of the solution. Only the strong consistency of these
estimators was proved. A comparison of different estimators of the
Hurst index for linear SDE driven by fBm was carried out in
\cite{kubmel2}. A rate of convergence of these estimators with
probability $1$ to the true value of a parameter was obtained in
\cite{kubmish}, when the maximum  length of the interval partition
tends to zero.

In this paper we will use a different approach to solving the
equation (\ref{e:diflygt}) than suggested by N.~T. Dung. We apply a
pathwise approach to stochastic integral equations and use a chain
rule for the composition of a smooth function and a function of
bounded $p$-variation with $1<p<2$. This approach allows us to simplify
the proof of existence and uniqueness of the solution of the
equation (1).  The goal of the paper is to construct strongly
consistent and asymptotically normal estimators of the Hurst parameter
and the volatility of SDE (\ref{e:diflygt}) using discrete observations of
the sample paths of the solution of the SDE. Estimators are constructed using
quadratic variations of second order increments. Estimators'
properties are determined asymptotically from the properties of
quadratic variations.

The paper is organized in the following way. In Section
\ref{s:main_result} we present the main results of the paper.
Section \ref{s:preliminaries} is devoted to several known results
needed for the proofs. In Section \ref{s:solution} we prove the
uniqueness of SDE (\ref{e:diflygt}) and give its explicit solution.
Sections \ref{s:proofs1}--\ref{s:proofs2} contain the proofs of
the main results. Finally, in Section \ref{s:modeling} some simulations are given in order to illustrate the obtained results.

\section{Main results}\label{s:main_result}

Let $X$ be a stochastic process satisfying (\ref{e:diflygt}). Denote
\[
\D^{(2)}_{n,k} X=X_{n,k+1}-2X_{n,k} +X_{n,k-1},\qquad \D_{n,k+1}
X=X_{n,k+1}-X_{n,k},
\]
where $X_{n,k}=X(t^n_k)$, $t^n_k=\frac{kT}{n}$. For the sake of
simplicity, we shall omit the index $n$ for points $t_k^n$ and for
notions $X_{n,k}$,  $\D^{(2)}_{n,k}$, $\D_{n,k}$, where the sample
size is equal to $n$. In case of different sample sizes the indices
will be retained.

In the formulation of Theorem \ref{t:case1} as well as in the sequel
we make use of the symbols $O_r, o_r$. Here is a brief explanation. Let
$(a_n)$ be a sequence of real numbers. $Y_n=O_r(a_n)$,
$a_n\downarrow 0$, means that there exists a.s. finite r.v.
$\varsigma$ with the property $\vert Y_n\vert\le \varsigma\cdot
a_n$; $Y_n=o_r(a_n)$ means that $\vert Y_n\vert\le \varsigma\cdot
b_n$ with $b_n=o(a_n)$. Particulary $Y_n=o_r(1)$ denotes the
sequence $(Y_n)$ which tends to $0$ a.s. as $n\to \infty$.

\begin{thm}\label{t:case} Assume that $X$ is a solution of SDE (\ref{e:diflygt}) with $a\le 0$. Let $X_0>0$, $c\neq 0$, and $1/2<H<1$. Then
\begin{align}
\widehat H^{(1)}_n&\tend{} H \quad\mbox{a.s.,}\nonumber\\
2\sqrt{n}\,\ln\frac{n}{T}\,(\widehat H^{(1)}_n-H)&\tendsd
\mathcal{N}(0,\sigma^2)\label{l:normalumas}
\end{align}
with a known variance $\sigma^2=\sigma^2(H)$ defined in Subsection
\ref{s:fbm}, where for $n>T$
\begin{align*}
\widehat H^{(1)}_n=&\varphi_{n,T}^{-1}\Bigg(\frac{1}{n}\sum_{k=1}^{n-1}\left( \frac{\Delta^{(2)}_k X}{c X_k}\right)^2\Bigg),\\
\varphi_{n,T}(x)=&\Big(\frac{T}{n}\Big)^{2x}(4-2^{2x})\quad\text{
 and  }\quad \varphi_{n,T}^{-1}\ \mbox{is the inverse
function of}\ \varphi_{n,T},\quad x\in(0,1).
\end{align*}
\end{thm}

\noindent Theorem \ref{t:case} gives an estimator of $H$ under the
assumption that $c$ is known. This is not always the case. Therefore
we present another estimator which is suitable when $c$ is unknown.

\begin{thm}\label{t:case1}
Let assumptions of Theorem \ref{t:case} hold and
\begin{equation*}
    \widehat H^{(2)}_n=\frac{1}{2}-\frac{1}{2\ln2}
    \ln\left(\frac{\sum_{k=1}^{2n-1}\left( \frac{\Delta^{(2)}_{2n,k}X}{ X_{2n,k}}\right)^2}{\sum_{k=1}^{n-1}\left( \frac{\Delta^{(2)}_{n,k}X}{
    X_{n,k}}\right)^2}\right).
\end{equation*}
Then
\begin{align}
\widehat H^{(2)}_n&\tend{} H \quad\mbox{a.s.,}\nonumber\\
2\ln 2\,\sqrt{n}(\widehat H^{(2)}_n-H)&\tendsd
\mathcal{N}(0,\sigma_*^2)
\end{align}
with a known variance $\sigma_*^2=\sigma_*^2(H)$ defined in
Subsection \ref{s:fbm}.
\end{thm}

\noindent To estimate $c$ one can use the last theorem provided
below.

\begin{thm}\label{t:case2} Assume that $1/2<H<1$ and $\widehat H^{(3)}_n=H+o_{r}(\phi(n))$ a.s. If $\phi(n)=o\left(\frac{1}{\ln n}\right)$ then
\begin{align*}
\widehat c^{\,2}_n=&\frac{n^{2 \widehat H^{(3)}_n-1}}{T^{2 \widehat
H^{(3)}_n}(4-2^{2 \widehat H^{(3)}_n})}\sum_{k=1}^{n-1}\left(
\frac{\Delta^{(2)}_kX}{ X_k}\right)^2\tend{} c^2 \quad\mbox{a.s.;}
\end{align*}
if $\phi(n)=o\left(\frac{1}{\sqrt{n}\ln n}\right)$ then
$\sqrt{n}\,(\widehat c^{\,2}_n-c^2)\tendsd \mathcal{N}(0,c^4\s^2)$.
\end{thm}
Hence estimation of $c$ and construction of confidence intervals is
always possible. However confidence intervals will shrink at the
rate which is a bit worse than the usual $n^{-1/2}$.

\section{Preliminaries}\label{s:preliminaries}

\subsection{Variation}\label{s:variation}

In this subsection we remind several facts about $p$--variation. For
details we refer the reader to \cite{DudleyNorvaisa-10}. Recall that
a  function $f:[a,b]\to\mathbb{R}$ has $p$--variation defined by
\begin{gather*}
    v_{p}(f;[a,b])=\sup_{\varkappa}\sum_{k=1}^{n}\abs{f(x_k)-f(x_{k-1})}^{p},
\end{gather*}
where $p\in(0,\infty)$ is a fixed number and supremum is taken over
all possible partitions
\[
    \varkappa:a=x_0<x_1<\dots<x_n=b,\qquad n\ge 1.
\]
If $v_p(f;[a,b])<\infty,$ $f$ is said to have a bounded
$p$-variation on $[a,b].$ Denote by $\mathcal{W}_p([a,b])$ the class
of functions defined on the $[a,b]$ with bounded $p$--variation. Let
$C\mathcal{W}_p([a,b])$ denote the class of functions that have
bounded $p$--variation and are continuous. Define
$V_p(f;[a,b])=v_p^{1/p}(f;[a,b])$, which is a seminorm on
$\mathcal{W}_p([a,b])$ provided $p\gs 1$ and $V_p(f;[a,b])$ is $0$
if and only if $f$ is a constant.

If $f\in \mathcal{W}_p([a,b])$, then f is bounded and $f\in
\mathcal{W}_{p_1}([a,b])$ for every $p_1>p\gs 1$. Let $0<p<\infty$
and $f,g\in \mathcal{W}_p([a,b])$. Then $fg\in
\mathcal{W}_p([a,b])$.

Let $f\in \mathcal{W}_q([a,b])$ and $h\in \mathcal{W}_p([a,b])$ with
$0<p<\infty$, $q>0,$ $1/p+\allowbreak 1/q>1.$ Then the integral
$\int_a^b f\,\mathrm{d}h$ exists as the Riemann--Stieltjes  integral
if $f$ and $h$ have no common discontinuities. If the integral
exists, the Love--Young inequality
\begin{equation}\label{in:LoveYoung}
\Bigg\vert \int\limits_a^bf\,\dr h-f(y)\big[ h(b)-h(a)
\big]\Bigg\vert \ls C_{p,q} V_q\big(f;[a,b]\big)V_p\big(h;[a,b]\big)
\end{equation}
holds for all $y\in [a,b]$, where $C_{p,q}=\zeta(p^{-1}+q^{-1})$ and
$\zeta(s)=\sum_{n\gs 1} n^{-s}.$

Let $f\in {\cal W}_q([a,b])$ and $h\in {\cal W}_p([a,b]).$ From
(\ref{in:LoveYoung}) it follows that
$$
V_p\Bigg (\int\limits_a^{\bt} f \,\dr h; [a, b] \Bigg) \leq C_{p, q}
V_{q,\infty}\big(f;[a,b]\big) V_p\big(h;[a, b]\big),
$$
where $V_{q,\infty}(f;[a,b])=V_q(f;[a,b])+\sup_{a\ls x\ls b}\vert
f(x)\vert$. Note that $V_{q,\infty}(f;[a,b])$ is a norm on ${\cal
W}_q([a,b])$, $q\gs 1$.

If $h\in C{\cal W}_p([a,b])$, then the indefinite integral $\int_a^y
f\,\dr h$, $y\in [a,b]$, is a continuous function.

To verify that some process solves the equation (\ref{e:diflygt}) we
use a chain rule for the composition of a smooth function and a
function of bounded $p$-variation with $0<p<2$. The chain rule
is based on the Riemann--Stieltjes integrals. As
usual, the composition $g\circ h$ on $[a, b]$ of two functions $g$ and $h$
is defined by $(g\circ h)(x)=g(h(x))$ whenever $h$ is defined on $[a, b]$ and
$g$ on the range of $h$.

\begin{thm}[Chain rule (see \cite{DudleyNorvaisa-10})]\label{t:salopek} Let $f = (f_1,\ldots, f_d): [a, b]\to \mathbb{R}^d$ be a function such that, for each $k = 1,\ldots,d$, $f_k\in C \mathcal{W}_p([a,b])$ with  $1\ls p < 2$. Let $g: \mathbb{R}^d \to \mathbb{R}$ be a differentiable function with locally Lipschitz partial derivatives $g^\prime_k$, $k=1,\ldots,d$. Then each $g^\prime_l\circ f$
is Riemann-Stieltjes integrable with respect to $f_k$ and
\[
(g\circ f)(b)-(g\circ f)(a)=\sum_{k=1}^d \int_a^b
(g_{k}^\prime\circ f)\, df_k.
\]
\end{thm}

\begin{prop}[Substitution rule (see \cite{DudleyNorvaisa-10})]\label{p: substitution rule} Let $f$, $g$, and $h$ be functions in $C \mathcal{W}_p([a,b])$, $1\ls p < 2$. Then
\[
\int_a^b f(x)\,d\bigg(\int_a^x g(y)\,dh(y)\bigg)=\int_a^b
f(x)g(x)\,dh(x).
\]
\end{prop}

\subsection{Several results on fBm}\label{s:fbm}

Recall that fBm $B^H=\{B^H_t\dvit t\gs 0\}$ with Hurst index $H\in
(0,1)$ is a real-valued continuous centered Gaussian process with
covariance given by \[ \E(B^H_t B^H_s)
=\frac12\big(s^{2H}+t^{2H}-|t-s|^{2H}\big).
\]
For consideration of strong consistency and asymptotic normality of
the given estimators we need several facts about $B^H$.

\medskip\noindent\emph{\textbf{Limit results.}} Let
\[
V_{n,T}=\frac{n^{2H-1}}{4-2^{2H}}\sum_{k=1}^{n-1}\big(T^{-H}\Delta^{(2)}_{n,k}
B^H\big)^2,\qquad H\neq \frac 12\,.
\]
Then (see \cite{Cohen-98}, \cite{IL}, \cite{begyn2})

\begin{gather*}
    V_{n,T}\tends{n\to\infty} 1\quad\text{ a.s.;} \\
    \sqrt{n}\left(
      \begin{array}{c}
        V_{n,T}-1 \\
        V_{2n,T}-1 \\
      \end{array}
    \right)\tendsd\mathcal{N}\left(\left(
                                     \begin{array}{c}
                                       0 \\
                                       0 \\
                                     \end{array}
                                   \right),
                                   \left(
               \begin{array}{cc}
                 \sigma^2(H) & \sigma_{1,2}(H) \\
                 \sigma_{1,2}(H) & \sigma^2(H)/2 \\
               \end{array}
             \right)
    \right)
\end{gather*}
with
\[
\sigma^2(H)= 2 + c_2(H) +
    c_1(H)\sum_{l=2}^{\infty}\big(\rho_{2-2H}(l)\big)^2,\qquad  \sigma_{1,2}(H)=2^{-2H}\big(3\sigma^2+\sigma_1^2+4\sigma_2^2\big),
\]
where
\begin{align*}
 c_1(H)=&\left(\frac{2H(2H-1)(2H-2)(2H-3)}{4-2^{2H}}\right)^2,\,\quad
    c_2(H)=\left(\frac{2^{2H+2}-7-3^{2H}}{4-2^{2H}}\right)^2,\\
    \rho_{2-2H}(l)=&\frac{\abs{l-2}^{2H}-4\abs{l-1}^{2H}+6\abs{l}^{2H}-4\abs{l+1}^{2H} +\abs{l+2}^{2H}}{2H(1-2H)(2-2H)(3-2H)}\,,\quad
    l\in\mathbb{Z},\\
     \sigma_1^2=&\frac{c_2(H)}{2} +
    c_1(H)\sum_{l=2}^{\infty}\rho_{2-2H}(l)\,\rho_{2-2H}(l-2),\\
    \sigma_2^2=&2 \sqrt{c_2(H)} +
    c_1(H)\sum_{l=2}^{\infty}\rho_{2-2H}(l)\,\rho_{2-2H}(l-1).
\end{align*}
These relationships imply that
\begin{gather*}
    \widetilde{H}_n=\frac{1}{2}-\frac{1}{2\ln
    2}\ln\left(\frac{V_{2n,T}}{2^{2H-1}V_{n,T}}\right)\to H\text{
    a.s.,}\\
    2\ln 2\,\sqrt{n}(\widetilde{H}_n-H)\tendsd\mathcal{N}(0;\sigma_*^2(H)),\qquad \sigma_*^2(H)=\frac 32 \sigma^2(H)-2\sigma_{1,2}(H).
\end{gather*}

\medskip\noindent\emph{\textbf{Variation of $B^H$.}} It is known that almost all sample paths of an fBm $B^H$ are locally H\"older of order strictly less than $H$, $0<H<1$. To be more
precise, for all ƒ$0<\eps<H$ and $T>0$, there exists a nonnegative
random variable $G_{\eps,T}$ such that $\mathbb{E}(\vert
G_{\eps,T}\vert^p)<\infty$ for all $p\gs 1$, and
\begin{equation}\label{in:fBmpokyt}
\vert B^H_t -B^H_s\vert \ls G_{\eps,T}\vert t-s\vert^{H-\eps}\qquad
a.s.
\end{equation}
for all $s,t\in [0,T]$. Thus $B^H\in C\mathcal{W}_{H_\eps}([0,T])$,
${H_\eps}=\frac{1}{H-\eps}$.

\medskip\noindent\emph{\textbf{SDE driven by fBm $B^H$.}}  A solution of the stochastic integral equation (\ref{e:diflygt}), on a given filtered probability space  $(\Omega,\,{\mathcal F},\,{\bf P},\,\mathbb{F})$ and with respect to the fixed fBm $(B^H,\,\mathbb{F}),$ $1/2<H<1,$ and
initial condition $\xi$, is an adapted to the filtration $\mathbb{F}$
continuous process $X=\{X_t\colon\ 0\leq  t\leq  T\}$ such that
$X_0=\xi$ a.s., ${\bf P}(\int_0^t\vert
f(X_s)\vert\,ds+\big\vert\int_0^t g(X_s)\,dB^H_s\big\vert<\infty)=1$
for every $0\leq  t\leq  T,$ and its almost all sample paths satisfy
(\ref{e:diflygt}).

\section{Solution of SDE  with polynomial drift}\label{s:solution}

First, we consider a non-random integral equation
\begin{equation}\label{e:diflygt1}
x_t=x_0+\int_0^t(a x_s^m+b x_s)\, ds + c \int_0^t x_s\, d h_s,\qquad
x_0\gs0,\quad a\ls 0,\ t\in[0,T],
\end{equation}
where $m\in\mathbb{N}$, $m\gs 2$, $h\in C\mathcal{W}_p([0,T])$, $1< p<2$.

\subsection{Uniqueness of the solution of equation (\ref{e:diflygt1})}

\begin{thm} The integral equation (\ref{e:diflygt1}) has a unique solution in the class of functions $C\mathcal{W}_p([0,T])$, $1< p<2$.
\end{thm}
\proof It is well-known that
\[
a^m-b^m=(a-b) \sum_{k=0}^{m-1}a^k b^{m-1-k}.
\]
Let us assume that we have two solutions $x_t$ and $y_t$ of equation
(\ref{e:diflygt1}). Since the solutions belong to the class of
functions  $C\mathcal{W}_p([0,T])$, $1< p<2$, then
\[
\vert x_t^m-y_t^m\vert\ls m L_{x,y,m,T}\vert x_t-y_t\vert,
\]
where
\[
L_{x,y,m,T}=\Big(\max\Big\{\max_{0\ls t\ls T}\vert x_t\vert,
\max_{0\ls t\ls T}\vert y_t\vert\Big\}\Big)^{m-1}.
\]
Further, one can find a subdivision
$\tau=\{0<\tau_1<\tau_2<\cdots<\tau_n=T\}$ such that
\[
V_p(h;[\tau_{k-1},\tau_k])\ls \big(4 \vert c\vert C_{p,p} \big)^{-1}
\]
for all $k$. Assume we have proved that
$x_{\tau_{k-1}}=y_{\tau_{k-1}}$. Then
\begin{align*}
V_{p,\infty}(x-y;[\tau_{k-1},\tau_k]) =&V_{p,\infty}(x-y-(x_{\tau_{k-1}}-y_{\tau_{k-1}});[\tau_{k-1},\tau_k])\\
\ls & 2\vert a\vert\int_{\tau_{k-1}}^{\tau_k} \vert
x^m_t-y^m_t\vert\,dt+2\vert b\vert\int_{\tau_{k-1}}^{\tau_k} \vert
x_t-y_t\vert\, dt\\
&+2\vert c\vert C_{p,p}
V_{p,\infty}(x-y;[\tau_{k-1},\tau_k])V_p(h;[\tau_{k-1},\tau_k])\\
\ls&\big(2L_{x,y,m,T} \vert a\vert+2\vert
b\vert\big)\int_{\tau_{k-1}}^{\tau_k} \vert x_t-y_t\vert\, dt\\
&+2\vert c\vert C_{p,p} V_{p,\infty}(x-y;[\tau_{k-1},\tau_k])
V_p(h;[\tau_{k-1},\tau_k])
\end{align*}
and
\begin{align*}
V_{p,\infty}(x-y;[\tau_{k-1},\tau_k])\ls& 4 \big(L_{x,y,m,T} \vert a\vert+\vert b\vert\big)\int_{\tau_{k-1}}^{\tau_k} \vert x_t-y_t\vert\, dt \\
\ls& 4 \big(L_{x,y,m,T} \vert a\vert+\vert
b\vert\big)\int_{\tau_{k-1}}^{\tau_k}
V_{p,\infty}(x-y;[\tau_{k-1},t])\, dt.
\end{align*}
By Gronwall's inequality we then have
$V_{p,\infty}(x-y;[\tau_{k-1},\tau_k])=0$.  Thus
$x_t=y_t$ on $[\tau_{k-1},\tau_k]$. Since
$x_{\tau_0}=x_0=y_{\tau_0}$ the claim of the theorem follows from
repetitive application of the above reasoning.

\subsection{Explicit solution of equation (\ref{e:diflygt1})}

\begin{thm} Function
\begin{equation}\label{e:solution}
x_t= e^{b t + c h_t}\bigg(x_0^{1-m}+(1-m)a \int_0^t e^{(m-1)(b s + c
h_s)}ds\bigg)^{1/(1-m)},\qquad t\in[0,T],
\end{equation}
is an element of the class of functions $C\mathcal{W}_p([0,T])$, $1<
p<2$, and  satisfies equation (\ref{e:diflygt1}).
\end{thm}
\proof We first show that $x\in
C\mathcal{W}_p([0,T])$, $1< p<2$. Let us denote $z_t=\exp\{b t+ c
h_t\}$ and $z_k=z_{t_k}$, where $\{0=t_0<t_1<\cdots<t_n=T\}$ is a
subdivision of the interval $[0,T]$. By the mean value theorem for the
function $e^x$, $x\in\mathbb{R}$, we obtain
\begin{align*}
\sum_{k=0}^{n-1}\vert z_{k+1}-z_k\vert^p\ls& \exp\big\{p \big[\vert b\vert T+2\vert c\vert V_p\big(h;[0,T]\big)\big]\big\}\sum_{k=0}^{n-1}\big\vert b\D_{k+1}+ c\D h_{k+1} \big\vert^p \\
\ls& 2\exp\big\{p \big[\vert b\vert T+2\vert c\vert
V_p\big(h;[0,T]\big)\big]\big\}\big[\vert b\vert^p T^p  +\vert
c\vert^p v_p\big(h;[0,T]\big)\big],
\end{align*}
where $\D_{k+1}=t_{k+1}-t_k$ and $\D h_{k+1}=h_{t_{k+1}}-h_{t_k}$.
Let
\[
f_t=x_0^{1-m}+(1-m)a \int_0^t e^{(m-1)(b s + c h_s)}ds.
\]
Repeatedly using the mean value theorem for the function
$u\mapsto u^{1/(1-m)}$, $u>0$, $m\gs 2$, we obtain
\begin{align*}
\sum_{k=0}^{n-1}\vert f^{1/(1-m)}_{k+1}-f^{1/(1-m)}_k\vert\ls&
x_0^{m}\sum_{k=0}^{n-1}\vert f_{k+1}-f_k\vert \ls
x_0^{m} v_1\big(f;[0,T]\big).
\end{align*}
Since  $z\in C\mathcal{W}_p([0,T])$ and $f\in
C\mathcal{W}_1([0,T])$, we have $x\in C\mathcal{W}_p([0,T])$.

Now we verify that function (\ref{e:solution}) satisfies the equation
(\ref{e:diflygt1}) by using the chain rule. Indeed, let $F(t,u,v)=e^{b t
+ c u}v^{1/(1-m)}$. Then by Theorem \ref{t:salopek} and Proposition
\ref{p: substitution rule} we get
\begin{align*}
x_t=&F(t,h_t,f_t)=F(0,h_0,f_0)+b\int_0^t x_s\,ds+c \int_0^t x_s\,dh_s\\
&+\frac{1}{1-m}\int_0^t e^{b s + c h_s}f_s^{1/(1-m)-1}\,d f_s\\
=&x_0+b\int_0^t x_s\,ds+c \int_0^t x_s\,dh_s+\frac{1}{1-m}\int_0^t e^{b s + c h_s}f_s^{m/(1-m)}\,df_s\\
=&x_0+b\int_0^t x_s\,ds+c \int_0^t x_s\,dh_s+\frac{1}{1-m}\int_0^t (1-m)a e^{m(b s + c h_s)}f_s^{m/(1-m)}\,ds\\
=&x_0+b\int_0^t x_s\,ds+c \int_0^t x_s\,dh_s+a\int_0^t x^m_s\,ds.
\end{align*}

\subsection{Solution of SDE}

Since almost all sample paths of $B^H$, $1/2<H<1$, are continuous
and have bounded
$H_\eps=\frac{1}{H-\eps}$-variation,
$1/2<H-\eps<1$, the pathwise Riemann-Stieltjes  integral
$\int_0^t X_s\, d B^H_s$ exists for $X\in
C\mathcal{W}_{H_\eps}([0,T])$. So SDE (\ref{e:diflygt}) makes sense
for almost all $\omega$. Thus the obtained result for a non-random
integral equation can be applied to an equation driven by a fBm.

\begin{thm} Suppose that $X_0>0$ and $m\gs 2$. The stochastic process
\begin{equation}\label{e:solution1}
X_t= e^{b t + c B^H_t}\bigg(X_0^{1-m}+(1-m)a \int_0^t e^{(m-1)(b s +
c B^H_s)}ds\bigg)^{1/(1-m)},\qquad t\in[0,T],
\end{equation}
for almost all $\omega$ belongs to  $C\mathcal{W}_{H_\eps}([0,T])$
and is the unique solution of (\ref{e:diflygt}).
\end{thm}

\section{Proof of the main Theorems}

\subsection{Proof of Theorem \ref{t:case}}\label{s:proofs1}

Before presenting the proof of this theorem, we  give an auxiliary
lemma. To avoid cumbersome expressions we extend notions of
$o_r,O_r$ to continuous time processes  as follows:
$Y_t=O_{r}(q(t))$ means that there exists some almost surely finite
non-negative r.v. $\zeta$ and $q(t)\tends{t\to 0+}
0+$, with the property $\abs{Y_t}\le \zeta q(t)$, $t\to 0$, whereas
$Y_t=o_{r}(q(t))$ means that $\abs{Y_t}\le \zeta w(t)$, $t\to 0$,
with $w(t)=o(q(t))$, $t\to 0$.

For $t\in[0,T)$ and $h\ge 0$ we also denote
\begin{equation*}
    \Delta Y_{t,t+h}=Y_{t+h}-Y_t,\qquad \Delta^{(2)} Y_{t,t+h}=\Delta Y_{t+h,t+2h}-\Delta
    Y_{t,t+h}.
\end{equation*}

\begin{lem}\label{l:nario_asimp} Assume that $H\in(1/2,1),a\le0,m\ge 2$ and $X$ is the solution of equation (\ref{e:diflygt}). Then for every $\eps>0$ such that $H>1/2+\eps$
\[
\Delta X_{t,t+h}=X_t\big(bh + c \Delta B^H_{t,t+h} +
O_r(h)\big),\qquad h\to 0+,\ t\in[0,T),
\]
and
\[
\Delta^{(2)}X_{t,t+h}=X_t\big(c \Delta^{(2)}B^H_{t,t+h} +
O_r(h^{2(H-\eps)})\big),\qquad h\to 0+,\ t\in[0,T).
\]
\end{lem}

\proof We rewrite the solution of equation (\ref{e:diflygt}) in the
following form
\[
X_t= X_0Z_t A_t,
\]
where
\[
Z_t=\exp\{b t+ c B^H_t\} \quad\mbox{and}\quad A_t =\bigg(1+(1-m)a
X_0^{m-1}\int_0^t Z_s^{m-1}ds\bigg)^{1/(1-m)}.
\]
Let $h\to 0+$. By (\ref{in:fBmpokyt}) and Maclaurin's expansion of
$e^{x}$
\begin{align}\label{e:Zt_asimptotika}
    Z_{t+h}=&e^{b(t+h)}e^{cB^H_{t+h}}=e^{bt}(1+bh+o(h^{2(H-\eps))})e^{cB^{H}_{t}+c\Delta B^H_{t,t+h}}\nonumber\\
    =&Z_t\left(1+bh+o(h^{2(H-\eps)})\right)\left(1+c\Delta B^H_{t,t+h}+ O_r(h^{2(H-\eps)})\right)\nonumber\\
    =&Z_t\left(1+ bh+ c\Delta B^H_{t,t+h} + O_{r}(h^{2(H-\eps)})\right).
\end{align}
Next note that $t\mapsto A_t\ge 1$ is non-decreasing whereas
\begin{equation}
\int_{t}^{t+h}Z_s^{m-1}ds=h Z_{t+\t h}^{m-1}=O_r(h),\qquad
\t\in(0;1),
\end{equation}
vanishes. Therefore Maclaurin's expansion of
$x\mapsto(1+x)^{\alpha}$ gives
\begin{align}\label{e:At_asimptotika}
    A_{t+h}=&\left(A_t^{1-m} +
    (1-m)aX_0^{m-1}\int_{t}^{t+h}Z_s^{m-1}ds\right)^{1/(1-m)}\nonumber\\
     =& A_t\left(1+I_{t,t+h}+o_r(h^{2(H-\eps)})\right),
\end{align}
where
\[
I_{t,t+h}=\frac{aX_0^{m-1}}{A_t^{1-m}}\int_{t}^{t+h}Z_s^{m-1}ds.
\]
Hence
\begin{align}\label{e:XT_asimptotika}
    X_{t+h}=&X_t\left(1+ bh+ c\Delta B^H_{t,t+h} +
    O_{r}(h^{2(H-\eps)})\right)\left(1+I_{t,t+h}+o_r(h^{2(H-\eps)})\right)\nonumber\\
    =& X_t\left(1+ bh+ c\Delta B^H_{t,t+h}  + I_{t,t+h}+
    O_{r}(h^{2(H-\eps)})\right).
\end{align}
Consequently,
\[
\Delta X_{t,t+h}=X_t\left(bh +c\Delta B^H_{t,t+h}
+I_{t,t+h}+O_{r}(h^{2(H-\eps)})\right)
\]
and taking into account (\ref{e:XT_asimptotika}) it follows that
\begin{align*}
    \Delta^{(2)}X_{t,t+h}=&X_{t+h}\left(bh+ c\Delta B^H_{t+h,t+2h} + I_{t+h,t+2h}+
    O_{r}(h^{2(H-\eps)})\right)\\
    &-X_t\big(bh +c\Delta B^H_{t,t+h} + I_{t,t+h}
    +O_{r}(h^{2(H-\eps)})\big)\\
    =&X_t\big(1+ bh+ c\Delta B^H_{t,t+h} + I_{t,t+h}
    +O_{r}(h^{2(H-\eps)})\big)\\
    &\times\big(bh+ c\Delta B^H_{t+h,t+2h} +
    I_{t+h,t+2h}+O_{r}(h^{2(H-\eps)})\big)\\
    &-X_t\left(bh+ c\Delta B^H_{t,t+h} + I_{t,t+h}
    +O_{r}(h^{2(H-\eps)})\right)\\
    =&X_t\left(c\Delta^{(2)}B^H_{t,t+h}+I_{t+h,t+2h}-I_{t,t+h}+O_{r}(h^{2(H-\eps)})\right).
\end{align*}

It remains to show that
$I_{t+h,t+2h}-I_{t,t+h}=O_{r}(h^{2(H-\eps)})$. To see this observe
that by (\ref{e:Zt_asimptotika})--(\ref{e:At_asimptotika}) and the mean
value theorem, for some $\theta_i\in(0;1),i=1,2$,
\begin{align*}
    I_{t+h,t+2h}-I_{t,t+h}=&\frac{aX_0^{m-1}}{A_{t+h}^{m-1}}\left(\int_{t+h}^{t+2h}Z_s^{m-1}ds-
    \left(\frac{A_{t+h}}{A_{t}}\right)^{m-1}\int_{t}^{t+h}Z_s^{m-1}ds\right)\\
    =&\frac{aX_0^{m-1}}{A_{t+h}^{m-1}}\left(\int_{t+h}^{t+2h}Z_s^{m-1}ds-
    \left(1+O_r(h)\right)^{m-1}\int_{t}^{t+h}Z_s^{m-1}ds\right)\\
    =&\frac{aX_0^{m-1}}{A_{t+h}^{m-1}}\left(\int_{t+h}^{t+2h}Z_s^{m-1}ds-\int_{t}^{t+h}Z_s^{m-1}ds+
    O_r(h^2)\right)\\
    =&\frac{aX_0^{m-1}}{A_{t+h}^{m-1}}\left(h Z_{t+h+\theta_1 h}^{m-1}-h Z_{t+\theta_2 h}^{m-1}+
    O_r(h^2)\right)\\
    =&\frac{aX_0^{m-1}}{A_{t+h}^{m-1}}\left(hZ_s^{m-1}\left(\left(\frac{Z_{s+\tau}}{Z_{\tau}}\right)^{m-1}-1\right)+
    O_r(h^2)\right)\\
    =&\frac{aX_0^{m-1}}{A_{t+h}^{m-1}}\left(hZ_s^{m-1}\left(\left(1+O_r(h^{H-\eps})\right)^{m-1}-1\right)+
    O_r(h^2)\right)\\
    =&\frac{aX_0^{m-1}}{A_{t+h}^{m-1}}\,O_r(h^{1+H-\eps})=O_r(h^{2(H-\eps)}),
\end{align*}
where $s=t+\theta_2h$ and $\tau=h(1+\theta_1-\theta_2)$.

\noindent\emph{\textbf{Proof of Theorem \ref{t:case}.}} Observe
first that the function
\[
\fy_{n,T}(x)=\Big(\frac{T}{n}\Big)^{2x}(4-2^{2x}),\qquad x\in(0,1),
\]
is continuous and strictly decreasing for $n>T$. Thus it has the
inverse function $\fy_{n,T}^{-1}$ for $n>T$. By Lemma
\ref{l:nario_asimp} we obtain
\begin{align*}
\frac{\fy_{n,T}(\widehat H^{(1)}_n)}{\fy_{n,T}(H)} =&\Big[\Big(\frac{T}{n}\Big)^{2H}(4-2^{2H})\Big]^{-1} \fy_{n,T}\Bigg( \fy_{n,T}^{-1}\Bigg(\frac{1}{n}\sum_{k=1}^{n-1}\left(\frac{\Delta^{(2)}_kX}{c X_k}\right)^2\Bigg)\Bigg)\nonumber\\ =&\Big[\Big(\frac{T}{n}\Big)^{2H}(4-2^{2H})\Big]^{-1}\Bigg(\frac{1}{n}\sum_{k=1}^{n-1} \left(\frac{\Delta^{(2)}_kX}{c X_k}\right)^2\Bigg)\nonumber\\
=&\frac{n^{2H-1}}{T^{2H}(4-2^{2H})}\bigg(\sum_{k=1}^{n-1} \left(\Delta^{(2)}_kB^H\right)^2 +O_r\Big(\frac{1}{n^{3(H-\eps)-1}}\Big)\bigg)\nonumber\\
=&\frac{n^{2H-1}}{4-2^{2H}}\sum_{k=1}^{n-1}
\left(T^{-H}\Delta^{(2)}_kB^H\right)^2+O_r\Big(\frac{1}{n^{H-3\eps}}\Big)
=V_{n,T}+O_r\Big(\frac{1}{n^{H-3\eps}}\Big)
\end{align*}
for $3\eps<H$.

Fix $\delta>0$. Let $\beta_n=\frac{\fy_{n,T}(\widehat
H^{(1)}_n)}{\fy_{n,T}(H)}$. Since $V_{n,T}\to 1$ a.s., the same
applies to $\beta_n$. Consequently, there exists $n_0(\omega)$ such
that
\begin{align*}
\fy_{n,T}(H+\delta)=&(4-2^{2(H+\delta)})\left(\frac{T}{n}\right)^{2(H+\delta)}<\beta_n
(4-2^{2H})\left(\frac{T}{n}\right)^{2H}\\
<&(4-2^{2(H-\delta)})\left(\frac{T}{n}\right)^{2(H-\delta)}=\fy_{n,T}(H-\delta)
\end{align*}
for all $n\ge n_0$, i.e. $\fy_{n,T}(\widehat
H^{(1)}_n)\in\big(\fy_{n,T}(H+\delta),\fy_{n,T}(H-\delta)\big)$ for
all $n\ge n_0$. Since the function $\fy_{n,T}$ is strictly
decreasing for all $n>T$ we conclude that $\widehat
H^{(1)}_n\in(H-\delta, H+\delta)$ for all $n\ge \max(n_0,T)$. Thus
$\widehat H^{(1)}_n$ is strongly consistent.

Now we prove the asymptotic normality of $\widehat H^{(1)}_n$. Note that
\begin{align}\label{e:stat_asimp}
\ln \bigg(\frac{\fy_{n,T}(\widehat H^{(1)}_n)}{\fy_{n,T}(H)}\bigg)= &\ln\Big(\frac{n}{T}\Big)^{-2(\widehat H^{(1)}_n-H)}+\ln\frac{4-2^{2\widehat H^{(1)}_n}}{4-2^{2H}}\nonumber\\
=&-2(\widehat H^{(1)}_n-H)\ln\Big(\frac{n}{T}\Big)
+\ln\big(4-2^{2\widehat H^{(1)}_n}\big)-\ln\big(4-2^{2H}\big).
\end{align}
By the Lagrange theorem for the function $h(x)=\ln(4-2^{2x})$,
$x\in(0,1)$, we obtain
\[
\ln\big(4-2^{2\widehat
H^{(1)}_n}\big)-\ln\big(4-2^{2H}\big)=(\widehat
H^{(1)}_n-H)h^\prime\big(H+\theta_n(\widehat H^{(1)}_n-H)\big)
\]
for some $\theta_n\in(0,1)$. Since $\widehat H^{(1)}_n\to H$ the
sequence $(h^\prime(H+\theta_n(\widehat H^{(1)}_n-H)))_{n\ge 1}$ is
bounded.

The equality (\ref{e:stat_asimp}) can be rewritten in the following way:
\begin{align}\label{l:isdestymas}
-2\sqrt{n}\,\ln\Big(\frac{n}{T}\Big)(\widehat H^{(1)}_n-H)=&
\frac{\sqrt{n}\Big(\ln\Big(\frac{\fy_{n,T}(\widehat
H^{(1)}_n)}{\fy_{n,T}(H)}\Big)-\ln
1\Big)}{\left(1-\frac{h^\prime\big(H+\theta_n(\widehat H^{(1)}_n-H)\big)}{2\ln\left(\frac{n}{T}\right)}\right)}\nonumber\\
=&\frac{\sqrt{n}\Big(\ln\Big(\frac{\fy_{n,T}(\widehat
H^{(1)}_n)}{\fy_{n,T}(H)}\Big)-\ln 1\Big)}{1+o_r(1)} \,.
\end{align}
Next note that the Delta method together with Slutsky's theorem and limit results of section \ref{s:fbm} imply
\begin{equation*}
\sqrt{n}\Big(\frac{\fy_{n,T}(\widehat
H^{(1)}_n)}{\fy_{n,T}(H)}-1\Big)
=\sqrt{n}\bigg(V_{n,T}-1\bigg)+O_r\Big(\frac{1}{n^{H-3\eps-1/2}}\Big)\\
\tendsd N(0;\sigma^2(H))
\end{equation*}
once $\eps>0$ is chosen to satisfy $3\eps<H-1/2$. It remains to
observe that assertion (\ref{l:normalumas}) immediately follows from
(\ref{l:isdestymas}) and repeated application of the Delta method along
with Slutsky's theorem.

\subsection{Proof of Theorem \ref{t:case1}}\label{s:proofs3}
Fix $\eps\in\left(0,\frac{H-1/2}{3}\right)$. By Lemma
\ref{l:nario_asimp}
\begin{align*}
\sum_{k=1}^{n-1}\left( \frac{\Delta^{(2)}_{n,k}X}{
X_{n,k}}\right)^2=&\sum_{k=1}^{n-1}\big(c \D^{(2)}_{n,k} B^H
+O_r\big(n^{-2(H-\eps)}\big)\big)^2\\
=&c^2 \sum_{k=1}^{n-1} \big(\D^{(2)}_{n,k} B^H\big)^2
+O_r\big(n^{1-3(H-\eps)}\big)\\
=&\frac{c^2T^{2H}(4-2^{2H})}{n^{2H-1}}\,V_{n,T}+O_r\big(n^{1-3(H-\eps)}\big)
\end{align*}
and in the same way
\begin{equation*}
    \sum_{k=1}^{2n-1}\left( \frac{\Delta^{(2)}_{2n,k}X}{ X_{2n,k}}\right)^2=
    \frac{c^2T^{2H}(4-2^{2H})}{(2n)^{2H-1}}\,V_{2n,T}+O_r\big(n^{1-3(H-\eps)}\big).
\end{equation*}
Since $V_{n,T},V_{2n,T}\to 1$ a.s. (see section \ref{s:fbm}),
\begin{align*}
    \widehat H^{(2)}_n=&\frac{1}{2}-\frac{1}{2\ln
    2}\ln\left(\frac{\frac{c^2T^{2H}(4-2^{2H})}{(2n)^{2H-1}}V_{2n,T} +O_r\big(n^{1-3(H-\eps)}\big)}{\frac{c^2T^{2H}(4-2^{2H})}{n^{2H-1}}V_{n,T} +O_r\big(n^{1-3(H-\eps)}\big)}\right)\\
    =&\frac{1}{2}-\frac{1}{2\ln
    2}\ln\left(\frac{V_{2n,T}}{2^{2H-1}V_{n,T}}\left(\frac{1+O_r\big(n^{-H+3\eps}\big)}{1 +O_r\big(n^{-H+3\eps}\big)}\right)\right)\\
    =&\widetilde H_n -\frac{1}{2\ln 2}
    \ln\left(1+O_r\big(n^{-H+3\eps}\big)\right)=\widetilde H_n+O_r\big(n^{-H+3\eps}\big).
\end{align*}
Now to finish the proof apply Slutsky's theorem and limit results of
section \ref{s:fbm}. Note that the limit variance $\sigma_*^2(H)$ of
$\widehat H^{(2)}_n$ equals to that of $\widetilde H_n$.

\subsection{Proof of Theorem \ref{t:case2}}\label{s:proofs2}

Let $0<\eps<1/2$ be such that $3\eps<H-1/2$ and
\[
\widetilde{c}^{\,2}_n=\frac{n^{2(H-\widehat
H^{(3)}_n)}(4-2^{2\widehat H^{(3)}_n)})}{T^{2(H-\widehat
H^{(3)}_n)}(4-2^{2H})}\,\widehat c^{\,2}_n.
\]
By Maclaurin's expansion and the mean value theorem we obtain
\begin{align*}
    \frac{n^{2(H-\widehat H^{(3)}_n)}(4-2^{2\widehat H^{(3)}_n)})}{T^{2(H-\widehat H^{(3)}_n)}(4-2^{2H})}=&
    \exp\left\{o_r(\phi(n))\ln\left(\frac{n}{T}\right)^2+\ln\big(4-2^{2\widehat H^{(3)}_n}\big)-\ln(4-2^{2{H}})\right\}\nonumber\\
    =&(1+o_r(\phi(n)\ln n))\exp\big\{o_r(\phi(n)) h^\prime\big(H+\theta_n(\widehat H^{(3)}_n-H)\big)\big\}\nonumber\\
        =&(1+o_r(\phi(n)\ln n))(1+o_r(\phi(n)))=(1+o_r(\phi(n)\ln n))
\end{align*}
for some $\theta_n\in(0,1)$, where $h(x)=\ln(4-2^{2x})$,
$x\in(0,1)$. Hence it suffices to show that $\widetilde{c}^{\,2}_n$
satisfies
\begin{equation*}
    \widetilde{c}^{\,2}_n\to c^2\quad\text{ a.s.,}\qquad \sqrt{n}(\widetilde{c}^{\,2}_n-c^2)\tendsd N(0;c^4\sigma^2)
\end{equation*}
and the claim will be proved. By Lemma \ref{l:nario_asimp}
\begin{align*}
\sum_{k=1}^{n-1}\left( \frac{\Delta^{(2)}_kX}{ X_k}\right)^2=c^2
\sum_{k=1}^{n-1} \big(\D^{(2)}_k B^H\big)^2
+O_r\big(n^{1-3(H-\eps)}\big).
\end{align*}
Thus
\begin{align*}
\widetilde{c}^{\,2}_n-c^2 =&c^2\bigg(\frac{n^{2H-1}}{T^{2H}(4-2^{2
H})}\sum_{k=1}^{n-1} \big(\D^{(2)} B^H_k\big)^2-1\bigg) +
O_r\big(n^{-H+3\eps}\big)\\
=& c^2\big(V_{n,T}-1\big)+ O_r\big(n^{-H+3\eps}\big).
\end{align*}
It remains to apply Slutsky's theorem and limit results stated in
section \ref{s:fbm}.

\section{Simulations}\label{s:modeling}
\setcounter{figure}{0}
The simulations of the obtained estimators presented below were performed using the R software environment \cite{R}. The performance of these estimators was assessed using sample paths of the Black-Scholes model, the Verhulst equation and the Landau-Ginzburg equation as defined below.
In Equation \ref{e:diflygt}, if we set $a=0$, $b=\l$, $c=\sigma$ we obtain the Black-Scholes model
\[
dX_t=\l X_t\, dt + \s X_t\, d B^H_t.
\]
Let us note that the results, obtained in Section \ref{s:main_result} are applicable for the Black-Scholes model. If we set $a=-1$, $b=\l$, $c=\s$, $m=2$ we obtain the Verhulst equation
\[
dX_t=\big(-X_t^2+\l X_t\big) dt + \sigma  X_t\, d B^H_t
\]
with the explicit solution
\[
X_t= \frac{\xi \exp\{\l t + \s B^H_t\}}{1+ \xi\int_0^t \exp\{\l s +
\s B^H_s\}ds}.
\]
If we set $a=-1$, $b=\l$, $c=\s$, $m=3$ we obtain the Landau-Ginzburg equation
\[
dX_t=\big(-X_t^3+\l X_t\big) dt + \sigma X_t\, d B^H_t.
\]
with the explicit solution
\[
X_t= \frac{\xi \exp\{\l t + \s B^H_t\}}{\sqrt{1+ 2\xi^2\int_0^t \exp\{2(\l s +
\s B^H_s)\}ds}}.
\]

Further on, unless explicitly stated otherwise, the (arbitrary) default values of $\xi = 3$, $b = 0.5$ and $c = 0.7$ were used, and in each case the conclusions were drawn from $500$ sample paths of the respective processes. The figures discussed further can be found in the Appendix.

Firstly in Figure \ref{fig:h_dep_h} we consider the dependance of numeric characteristics of the estimators $\widehat H^{(1)}_n$ (defined in Theorem \ref{t:case} and denoted in the figures as H1) and $\widehat H^{(2)}_n$ (defined in Theorem \ref{t:case1} and denoted in the figures as H2) on the true value of the parameter $H$.

\begin{figure}[h]
\centering
\includegraphics[scale=0.62]{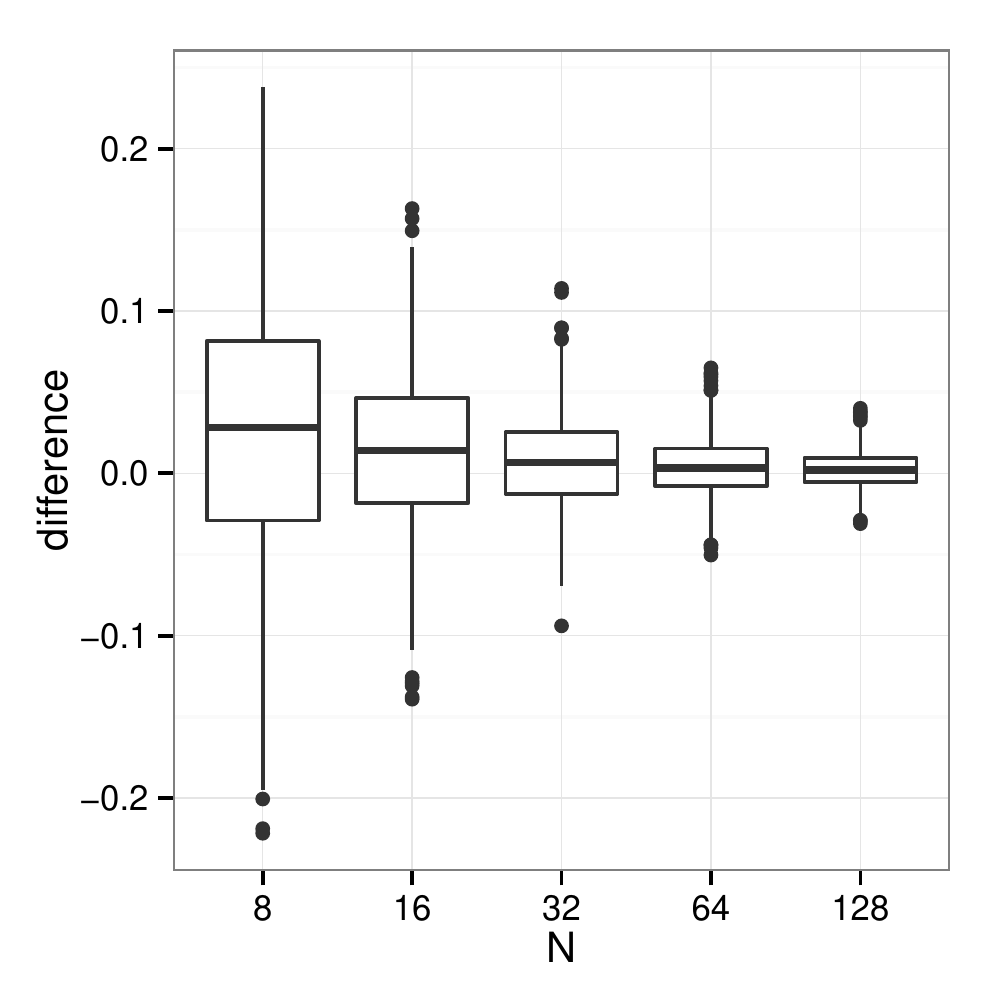}
\caption{Behavior of $\widehat H^{(1)}_n$ in case of short sample paths}
\label{fig:h1_small_n}
\end{figure}

The results imply that the accuracy of $\widehat H^{(1)}_n$ surpasses that of $\widehat H^{(2)}_n$ by at least an order of magnitude, and becomes even more precise for the values of $H$ close to 1. The practical use of $\widehat H^{(1)}_n$, however, requires the knowledge of the true value of the parameter $c$, which is not always available. However if that is the case then the estimator $\widehat H^{(1)}_n$ yields accurate estimators of $H$ even for very short sample paths, as illustrated in Figure \ref{fig:h1_small_n}. It is also worth to note that the estimators do not display notable dependance on the underlying process.

\begin{figure}[h]
\centering
\includegraphics[scale=0.7]{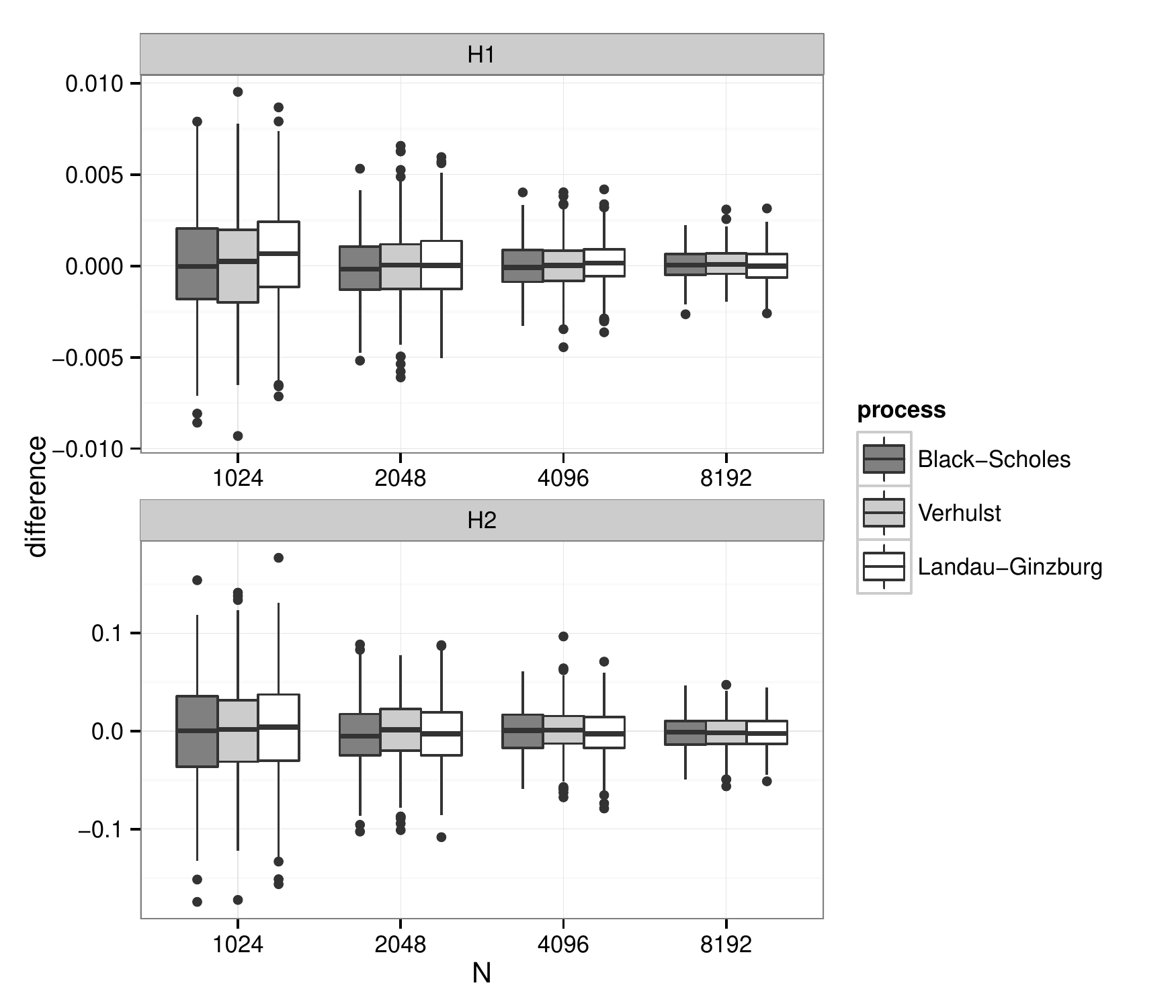}
\caption{The differences between the estimated and real values of $H$}
\label{fig:h_dep_n}
\end{figure}

The next question on the agenda is, how rapidly does the accuracy of the estimators grow as the length of sample paths is increased? Figure \ref{fig:h_dep_n} presents the boxplots of both estimators as the length of the sample paths, $N$, varies from 1024 to 8192 points. It can be seen that the estimator $\widehat H^{(1)}_n$ is considerably more precise than $\widehat H^{(2)}_n$, yet the rate at which their accuracy grows with increased sample sizes is roughly the same. For both of these estimators the inter-quartile ranges shrink by roughly 30\% as the sample length is doubled.

\begin{figure}[h]
\centering
\begin{minipage}{0.45\textwidth}
\centering
\includegraphics[scale=0.6]{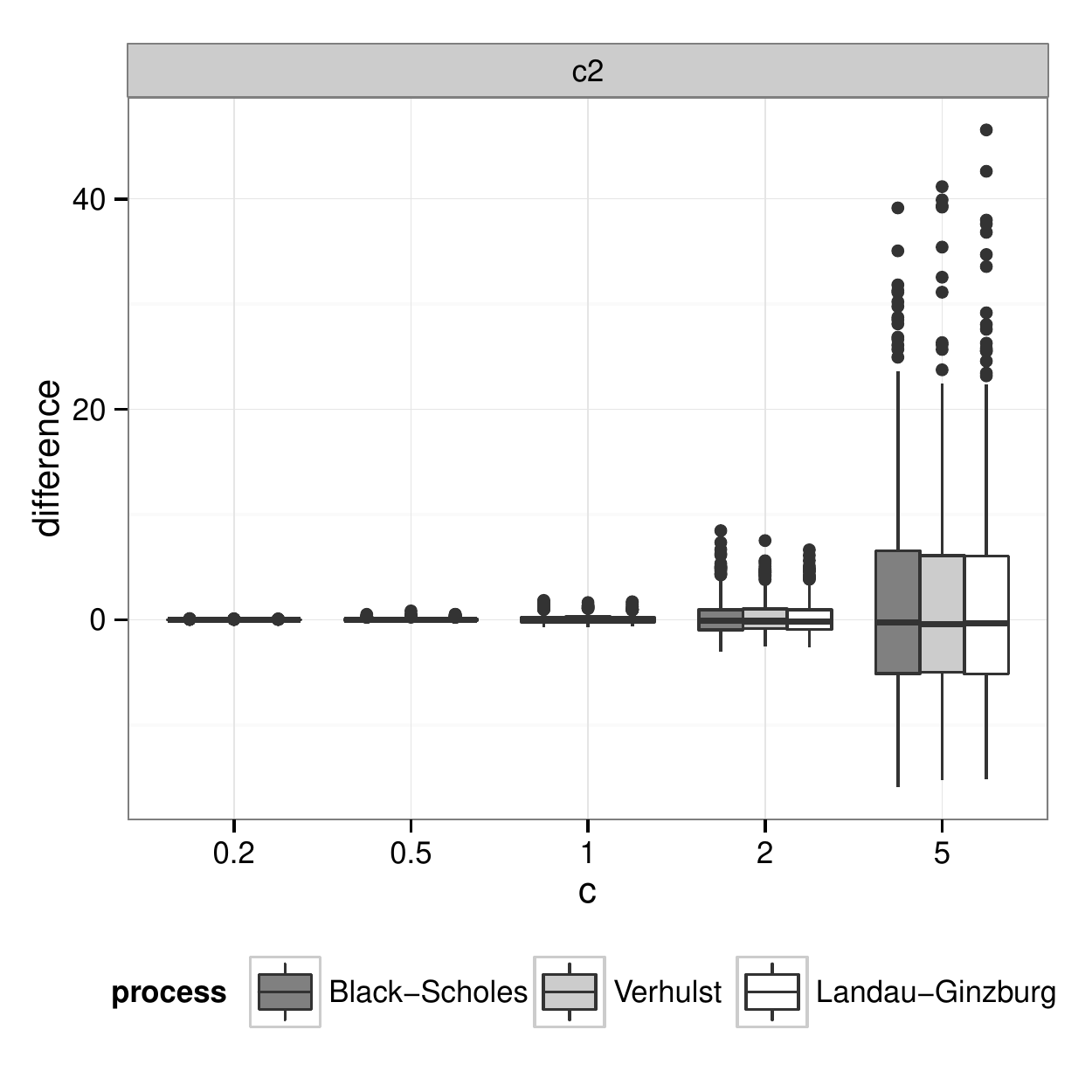}
\caption{$\widehat c^{\,2}_n - c^2$, dependance on $c$}
\label{fig:c_dep_c}
\end{minipage}\hfill
\begin{minipage}{0.45\textwidth}
\centering
\includegraphics[scale=0.7]{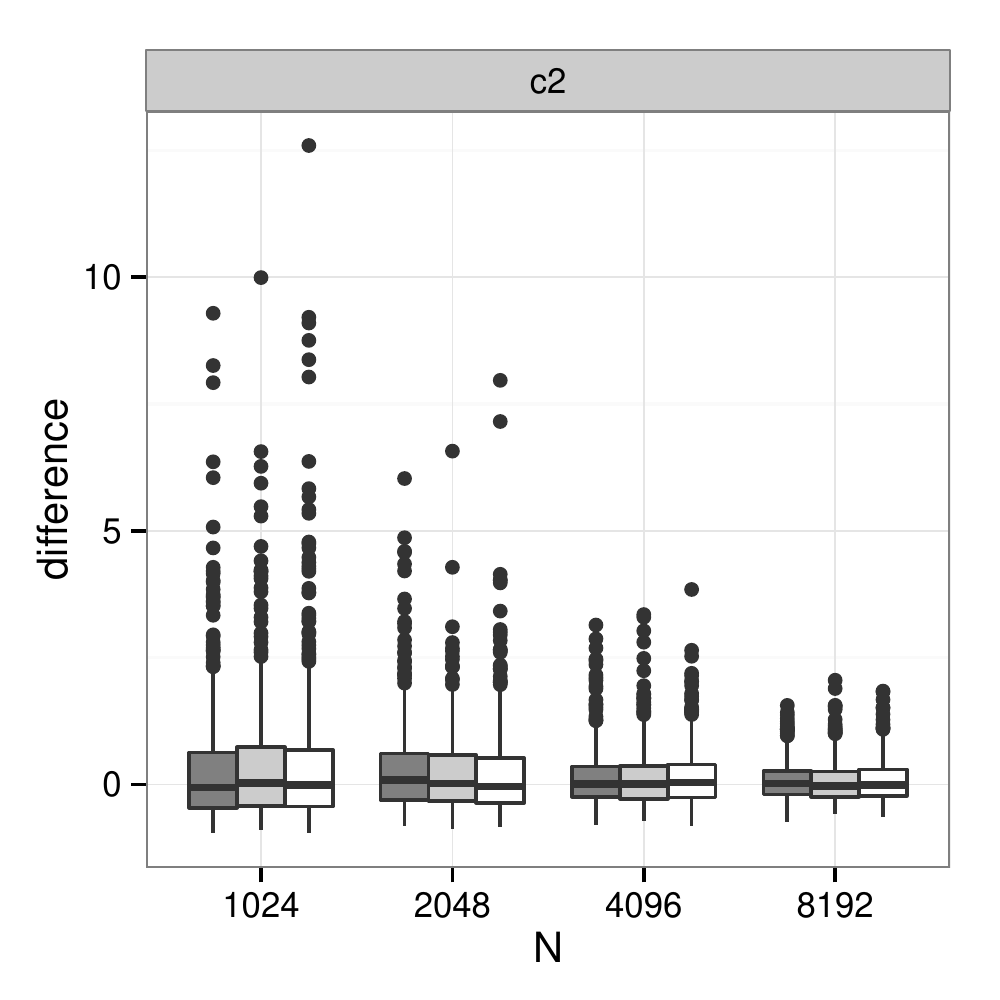}
\caption{$\widehat c^{\,2}_n - c^2$, dependance on $N$}
\label{fig:c_dep_n}
\end{minipage}
\end{figure}

Now let's examine the behavior of $\widehat c^{\,2}_n$, the estimator of the volatility defined in Theorem \ref{t:case2} and denoted c2 in the figures. As before, we present the boxplots of this estimator both for different values of the true parameter $c$ (Figure \ref{fig:c_dep_c}) and different sample path lengths (Figure \ref{fig:c_dep_n}). We can observe that the distribution of the estimator is somewhat right-skewed, which is hardly surprising since we're estimating $c^2$, not $c$ itself. Also, Figure \ref{fig:c_dep_c} shows that the variance of $\widehat c^{\,2}_n$ grows rapidly as the true value of the parameter increases. This is further illustrated by Table \ref{tab:c_dep_c} presenting the average bias $\widehat c^{\,2}_n - c^2$ along with the average variance of this estimator for different values of $c$. The adjusted R-squared of the linear model $Var(\widehat c^{\,2}_n) = k c^4 + b$ is $0.9838$ which is in agreement with the asymptotics presented in Theorem \ref{t:case2}.

\begin{table}[h]
\centering
\begin{tabular}{c|ccccc}
$c$ & 0.2 & 0.5 & 1 & 2 & 5 \\ \hline
Bias & 0.006 & 0.031 & 0.133 & 0.460 & 3.188 \\
Variance & 0.001 & 0.019 & 0.306 & 5.159 & 193.4
\end{tabular}
\caption{Numeric characteristics of $\widehat c^{\,2}_n$}
\label{tab:c_dep_c}
\end{table}

\newpage
\section*{Appendix: Figures}

\begin{figure}[h]
\centering
\includegraphics[scale=0.65]{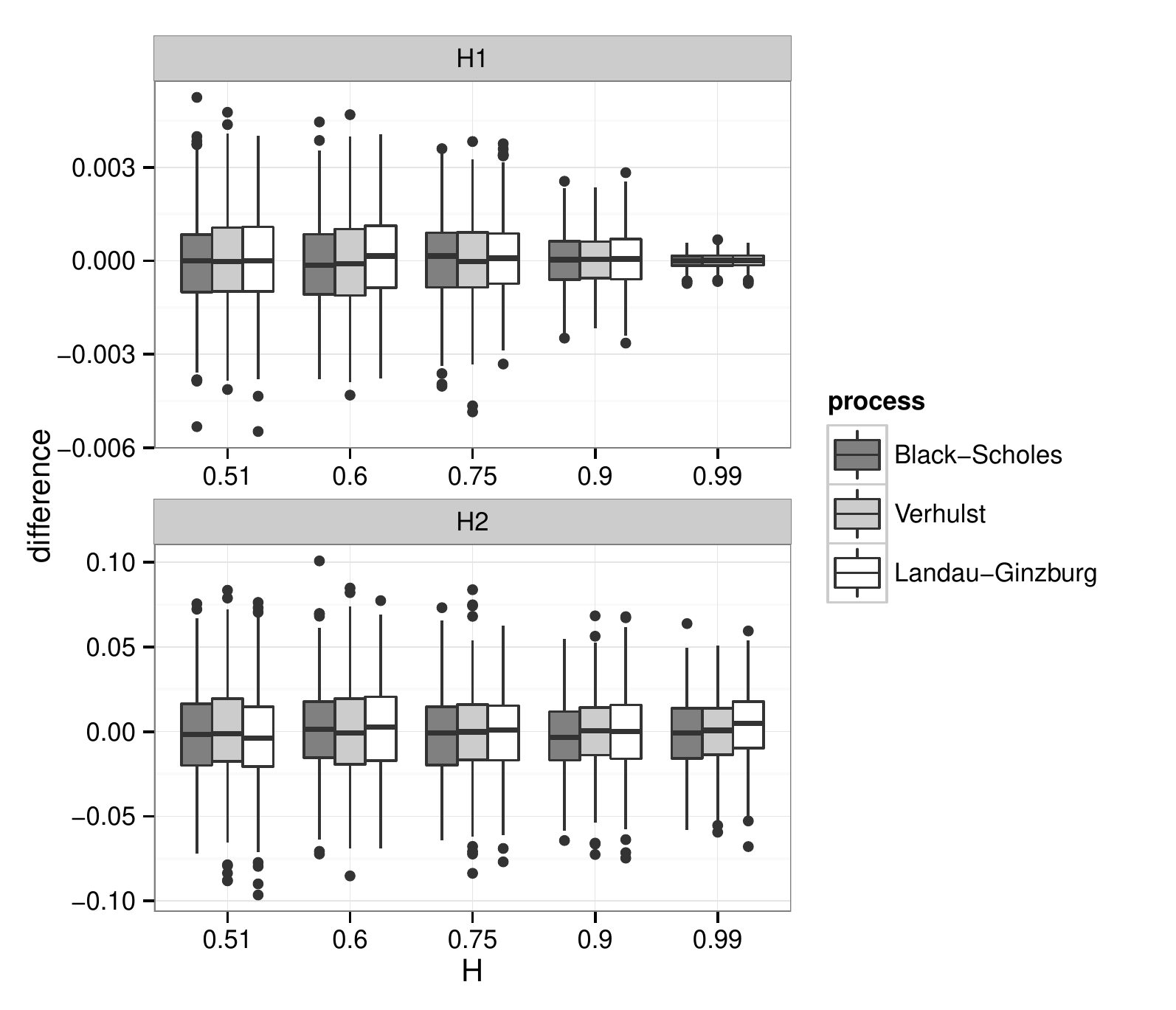}
\caption{The differences between the estimated and real values of $H$}
\label{fig:h_dep_h}
\end{figure}

\end{document}